\documentclass[10pt,reqno]{amsart}

\usepackage{amsmath,amsfonts,amssymb,pstricks,booktabs, enumerate,hyperref,multicol,graphicx,subfig,float,caption,subfig}

\begin{document}
\title{ Lie symmetry analysis for similarity reduction and solutions of $(3+1)$- dimensional Calogero-Bogoyavlenskii-Schiff Equation }
\author{}
\maketitle
\begin{center}
Vishakha Jadaun, Sachin Kumar
\end{center}

\begin{abstract}
 It is shown that the novel Lie group of transformations method is a competent and prominent tool in solving nonlinear partial differential equations(PDEs) in mathematical physics. Lie group analysis is used to carry out the similarity reduction and exact solutions of the $(3+1)$-dimensional Calogero-Bogoyavlenskii-Schiff (CBS) equation. This research deals with the similarity solutions of CBS equation. We have obtained the infinitesimal generators, commutator table of Lie algebra, symmetry group and similarity reduction for the CBS equation. For the different Lie algebra, Lie symmetry method reduced $(3+1)$-dimensional CBS equation into new $(2+1)$-dimensional partial differential equations and again using Lie symmetry method these PDEs are reduced into various ordinary differential equations(ODEs).

\end{abstract}
 \textbf{keywords:} $(3+1)$-dimensional Calogero-Bogoyavlenskii-Schiff equation, Lie symmetries, Similarity transformations method,Infinitesimal generator, Similarity solutions. 
\\

\section{Introduction}

Non-linear partial differential equations(PDEs) exhibit a rich variety of non-linear phenomena, arise in many physical fields like the stratified shear flow in ocean and atmosphere, condense matter physics, fluid mechanics, etc. Therefore, seeking exact solutions of  non-linear PDEs specifically the non-linear evolution equations (NLEEs) play an important role to look into the internal mechanism of convoluted physical phenomena. Most of the physical phenomena such as, fluid mechanics, quantum mechanics, electricity, plasma physics, chemical kinematics, propagation of shallow water waves and optical fibers are modelled by non-linear evolution equation and the appearance of solitary wave solutions in nature is somewhat frequent. The non-linear processes are one of the major challenges and not easy to control because the non-linear characteristic of the system abruptly changes due to some small changes in valid parameters including time.

 In the literature, many significant methods have been proposed for obtaining exact solutions of non-linear partial differential equations (PDEs) such as the Exp-function method, the Jacobi-elliptic method, the Lie B\"acklund transformations, the rational sine-cosine method, the Lie group of transformation method, the Hirota's method, Hirota bilinear forms, the tanh-sech method and so on \cite{GJ1,GJ7,OB10}. Lie group method of infinitesimal transformations \cite{ GPJM12, GS6, PO9, LVO8} which, has always been and still is, a great tool to find the analytical solution of non-linear partial differential equations(PDEs). A number of mathematicians have used this tool in many areas of scientific fields such as solid state physics, plasma physics, fluid dynamics, mathematical biology and chemical kinetics. S. Sahoo et.al \cite{SG} recently presented their work for modified Kdv-Zakharov-Kuznetsov equation by using Lie group method.
 
  In this paper we will study the generalized (3+1)-dimensional Calogero-Bogoyavlenskii-Schiff (CBS) equations
\begin{align*}
v_t+& \Phi(v)v_y+\Phi_1(v) v_z=0\\ \nonumber
\Phi(v)&= \partial^2_x +a v+b v_x \partial^{-1}_x\\ \nonumber
\Phi_1(v)&= \partial^2_x +c v+d v_x \partial^{-1}_x 
\end{align*}
or equivalently
\begin{equation}\label{c}
 v_t+a v v_y+c v v_z+b v_x \partial^{-1}_x v_y+d v_x  \partial^{-1}_x v_z+v_xxy+v_xxz=0
\end{equation}
where $a,b,c,$ and $d$ are parameters. Using a dimensional reduction $\partial z = \partial y =  \partial x$, (\ref{c}) will be reduced to the standard Korteweg-de Vries (KdV) equation. The
$(3+1)$-dimensional CBS equation (\ref{c}) can be obtained in the potential form  
 \begin{align}{\label{eq1}}
\Delta := u_{xt}+a u_x u_{xy}+b u_y u_{xx}+c u_x u_{xz}+d u_z u_{xx}+u_{xxxy}+u_{xxxz}=0
 \end{align}
 using the potential $v = u_x$. The CBS equation was first constructed by Bogoyavlenskii and Schiff in different ways. Bogoyavlenskii used the modified Lax formalism \cite{Y,TK,XYZ}, whereas Schiff derived the same equation by reducing the self-dual Yang-Mills equation \cite{TK, DW, XD,MM}.
 In many studies,  the forms of the arbitrary parameters which appear in the underlying model are assumed. However, the Lie symmetry approach through the method of group classification has proven to be a powerful tool in specifying the forms of these parameters naturally. The Lie group method is powerful technique to construct the exact solution of non-linear PDEs. Furthermore, based on the Lie group method, many types of exact solutions of PDEs can be considered, such as the traveling wave solutions, similarity solutions, soliton solutions, fundamental solutions, and so on.\\

\section{Method of symmetries}

In this section, we recall the general procedure for determining symmetries for any system of partial differential equation. Let us consider the general case of partial differential equations of order $p$ with $m-$dependent variables and $n-$independent variables is given as a system of equations
 \begin{equation}{\label{sys}}
\Delta_\textsl{s}(x,u^p)=0, \qquad \textsl{s}=1, \ldots,\emph{l}
\end{equation}
Here, $x=(x_1,x_2 \ldots x_n)$, $u=(u_1,u_2 \ldots u_m)$ and the derivatives of $u$ with respect to $x$ upto $p$, where $u^p$ express all the derivatives of $u$ of all orders from $0$ to $p$. we consider one parameter Lie group of infinitesimal transformation acting on dependent variable $u$ and independent variables $(x,y,t)$ of (\ref{sys}).
\begin{equation}
\left.
\begin{split}
\check{x^i} &=x^i+\varepsilon \xi^i(x,u)+O(\varepsilon^2), \quad i=1,2, \ldots, n\\
\check{u^j} &= u^j+\varepsilon \eta^j(x,u)+O(\varepsilon^2),  \quad j=1,2, \ldots, m
\end{split}
\right\}
\end{equation}
  where $\xi^i $ and $\eta^j$ are infinitesimals of transformation for the independent and dependent variables, respectevely and  $(\varepsilon)$ is  group  parameter which  is  admitted by  the  system  (\ref{sys})\\
  The infinitesimal generator $\textbf{v}$ associated with the above group of transformation can be written as
  \begin{equation}
  \textbf{v}= \sum_{i=1}^n \xi^i(x,u)\partial_x+\sum_{\alpha=1}^m \eta^\alpha(x,u) \partial_{u^\alpha}
  \end{equation}
Lie  group  of transformations are such  that  if $ u $ is  a  solution of system (\ref{sys}) then $ \check{u}$ is  also a solution. The method for finding group symmetry \cite{LVO8} is by finding corresponding infinitesimal generator of Lie group of  transformations. This yields to overdetermined, linear system of equation for infinitesimals $\xi^i(x,u), \eta^j(x,u)$. The invariance of system (\ref{sys}) under the infinitesimal transformation leads to the invariance conditions 
 \begin{equation}
 Pr^{(p)} \textbf{v}[\Delta_\textsl{s}(x,u^{(p)})] =0,\textsl{s}=1,\ldots,\emph{l} \quad \text{whenever} \quad \Delta_\textsl{s}(x,u^{(p)})=0
 \end{equation}
 where $Pr^{(p)}$ is called $p^{th}$-order prolongation vector field $\textbf{v}$ is given by
 \begin{equation}
 Pr^{(p)} \textbf{v}=\textbf{v}+\sum_{\alpha=1}^n \sum \eta_\alpha^\textit{j}(x,u^{(p)}) \partial u_\textit{j}^\alpha
\end{equation}
where $\textit{j}=(\textit{j}_1, \ldots, \textit{j}_k),  1\leq\textit{j}_k \leq n,  1\leq k \leq p$  , and sum is all over the orders of \textit{j}'s of order $0<\textit{j}\leq p$. If \textit{j}=k, the coefficients $\eta_\alpha^\textit{j}$ of $ \partial u_\textit{j}^\alpha$ will only depend on $k^{th}$ and lower order derivatives of $u$;
\begin{equation}
\eta_\alpha^\textit{j}(x,u^{(p)}= D_\textit{j}\left( \eta_\alpha - \sum_{\textit{i}=1}^n \xi^\textit{i} u_\textit{i}^\alpha\right)+\sum_{\textit{i}=1}^n  \xi^\textit{i}u_{\textit{j},\textit{i}}^\alpha.
\end{equation}
where $u_\textit{i}^\alpha= \frac{\partial u^\alpha}{\partial x^\textit{i}}$,  $u_{\textit{j},\textit{i}}^\alpha= \frac{\partial u_\textit{i}^\alpha}{\partial x^\textit{i}}$.\\
The set of all infinitesimal symmetries of this system have one important property which is that they form a Lie algebra under the usual Lie bracket.\\
\\

\section{ Lie symmetry analysis of $(3+1)-$dimensional CBS equation}

First we sketch the derivation of the Lie similarity reductions of the CBS equation using Lie symmetry analysis as given in \cite{GJ7}. we consider the one parameter Lie group of infinitesimal transformation on $(x_1=x, x_2=y, x_3=z, x_4=t, u_1=u)$,
\begin{align}
\check{x} &= x+\varepsilon \xi^1(x,y,z,t,u)+O(\varepsilon^2)\\
\check{y}&=y+\varepsilon \xi^2(x,y,z,t,u)+O(\varepsilon^2)\\
\check{z} &= t+\varepsilon \xi^3(x,y,z,t,u)+O(\varepsilon^2)\\
\check{t} &= t+\varepsilon \tau(x,y,z,t,u)+O(\varepsilon^2)\\
\check{u} &= u+\varepsilon \eta(x,y,z,t,u)+O(\varepsilon^2)
\end{align}
where $\varepsilon$ is the continuous group parameter and $\xi^1, \xi^2, \xi^3, \tau$ and $\eta$ are the infinitesimal of the transformations for the independent and dependent variable respectively, which are yet to be determined. The associated vector field of the form:\\
\begin{equation*}
\textbf{v}=\xi^1(x,y,z,t,u) \partial_x+\xi^2(x,y,z,t,u) \partial_y+\xi^3(x,y,z,t,u) \partial_z+\tau(x,y,z,t,u) \partial_t+\eta(x,y,z,t,u) \partial_u
\end{equation*}
\\
Using the invariance condition $Pr^{(4)}\textbf{v}(\Delta)=0$ whenever $\Delta=0$ and $Pr^{(4)}\textbf{v}$ is the fourth prolongation of $\textbf{v}$, thus the Infinitesimal criteria	 for the invariance  of Eq. \ref{eq1} would be:
\begin{equation}\label{a}
 \eta_{xt}+a \eta_x u_{xy}+a u_x \eta_{xy}+b\eta_y u_{xx}+b \eta_{xx} u_x+c \eta_x u_{xz}+c \eta_{xz}u_x+d \eta_z u_{xx}+ d \eta_{xx} u_z+\eta_{xxxy}+\eta_{xxxxz}
\end{equation}
Applying the fourth prolongation of $\textbf{v}$ to Eq.(\ref{a}), one can obtain a large overdetermined system of coupled partial differential equations which are called determining equations
\begin{align}\label{s}
\xi^1_u&=\xi^1_y=\xi^1_z=\xi^1_x-\frac{\xi^3_z}{2}-\frac{\tau_t}{2}=0 \nonumber \\
\xi^2_u&=\xi^2_z= \xi^2_y- \xi^2_z = \xi^2_t- \frac{a \xi^3_t}{c} =0 \nonumber \\
\xi^3_u& = \xi^3_x = \xi^3_y = \xi^3_{tt} = \xi^3_{tz}=0\nonumber \\
\tau_x& =\tau_y = \tau_z = \tau_u = \tau_{tt}=0 \nonumber \\
\eta_u &= \frac{-\tau_t}{2}+\frac{\xi^3_z}{2}, \quad \eta_x = \frac{\xi^3_t}{c}, \quad \eta_y = \frac{-d \eta_z +\xi^1_t}{b} 
\end{align} 
  The infinitesimals can be find by solving ``determining equations"(\ref{s})  which yields the following infinitesimals:
\begin{eqnarray}{\label{eq2}}
\begin{split}
\xi^1 =& (c_{1}-c_{4}) x + \lambda(t) \\
\xi^2 =&\quad 2 c_{4} y + \frac{a}{c}c_{3}t+c_{6} \\
\xi^3 =&\quad 2c_{3}t +2 c_{4} z+c_{5}\\
\tau =&\quad 2c_{1}t +c_{2}\\
\eta =& -(c_{1}-c_{4}) u + \frac{c_{3}}{c} x+ \frac{\lambda^\prime(t)}{b} y + \gamma\left(\frac{b z -d y}{b}, t\right)
\end{split}
\end{eqnarray}
where $c_{1}$, $c_{2}$, $c_{3}$, $c_{4}$, $c_{5}$ and $c_{5}$ are arbitrary constants.  $\gamma(t)$ and $\lambda(t)$ are arbitrary functions of t. The prime $(\prime)$ denotes the differentiation  with respect to it's indicated variable throughout the paper.\\
Hence, the infinitesimal generators of the corresponding Lie algebra are given by
\begin{eqnarray}{\label{eq3}}
\begin{split}
\textbf{v}_{1} =& \quad (x+\lambda(t))\frac{\partial}{\partial x}+ 2 t \frac{\partial}{\partial t}+\left(-u+\frac{\lambda^\prime(t) y}{b} +\gamma \left(\frac{b z-d y}{b},t\right)\right) \frac{\partial}{\partial u}\\
\textbf{v}_{2} =& \quad  \lambda(t) \frac{\partial}{\partial x}+\frac{\partial}{\partial t}+\left(\frac{\lambda^\prime(t) y}{b} +\gamma \left(\frac{b z-d y}{b},t\right)\right)\frac{\partial}{\partial u}\\
\textbf{v}_{3} =& \quad \lambda(t) \frac{\partial}{\partial x}+\frac{a t}{c} \frac{\partial}{\partial y}+t \frac{\partial}{\partial z}+\left(\frac{x}{c}+\frac{\lambda^\prime(t) y}{b} +\gamma \left(\frac{b z-d y}{b},t\right)\right) \frac{\partial}{\partial u}\\
\textbf{v}_{4} =& \quad(-x+\lambda(t))\frac{\partial}{\partial x}+ 2 y \frac{\partial}{\partial y}+2 z \frac{\partial}{\partial z}\left(u+\frac{\lambda^\prime(t) y}{b} +\gamma \left(\frac{b z-d y}{b},t\right)\right) \frac{\partial}{\partial u}\\
\textbf{v}_{5} =& \quad \lambda(t) \frac{\partial}{\partial x}+\frac{\partial}{\partial z}+\left(\frac{\lambda^\prime(t) y}{b} +\gamma \left(\frac{b z-d y}{b},t\right)\right)\frac{\partial}{\partial u}\\
\textbf{v}_{6} =& \quad \lambda(t) \frac{\partial}{\partial x}+\frac{\partial}{\partial y}+\left(\frac{\lambda^\prime(t) y}{b} +\gamma \left(\frac{b z-d y}{b},t\right)\right)\frac{\partial}{\partial u}
\end{split}
\end{eqnarray}
It is convenient to display the commutators of a Lie algebra through its commutator table whose $(i,j)^{th} $ entry is $[\textbf{v}_{i}, \textbf{v}_{j}]=\textbf{v}_{i} \ast \textbf{v}_{j}=\textbf{v}_{i} \cdot  \textbf{v}_{j}-\textbf{v}_{j} \cdot\textbf{v}_{i}$ . The commutator table is antisymmetric with its diagonal elements all zero as we have $[\textbf{v}_{\alpha}, \textbf{v}_{\beta}]= -[\textbf{v}_{\beta}, \textbf{v}_{\alpha}]$ (for more details see \cite{GJ7, GS6}). The structure constants are easily read off from the commutator table.

For the infinitesimal generators (\ref{eq3}) we have the following commutator table: 
\begin{table}[H]
\centering
\begin{tabular}{c|cccccc}
$\ast$ & $\textbf{v}_1$ & $\textbf{v}_2$ & $\textbf{v}_3$ & $\textbf{v}_4$ & $\textbf{v}_5$& $\textbf{v}_6$ \\
\hline
$\textbf{v}_1$&0 & $-2\textbf{v}_2$ & $2 \textbf{v}_3$ & 0 & 0 & 0\\
$\textbf{v}_2$ & $2\textbf{v}_2$ &0& $\frac{a}{c} \textbf{v}_6+\textbf{v}_5$ & 0 & 0 & 0\\
$\textbf{v}_3$ &$-2\textbf{v}_3$ & $-\frac{a}{c} \textbf{v}_6-\textbf{v}_5$ & 0 &  $2 \textbf{v}_3$& 0&0\\
$\textbf{v}_4$ &0&0&$-2 \textbf{v}_3$ &0&$-2 \textbf{v}_5$&$-2 \textbf{v}_6$\\
$\textbf{v}_5$ &0&0&0&$ 2 \textbf{v}_5$ &0&0\\
$\textbf{v}_6$ &0&0&0&$2 \textbf{v}_6$&0&0\\
\hline
\end{tabular}
\end{table}
Here, it is clear that the CBS equation contain infinite continuous group of transformations which is generated by the infinite-dimensional Lie algebra spanned by vector fields (\ref{eq3}). These generators are linearly independent. In general, there are an infinite number of subalgebras for this Lie algebra formed from linear combinations of generators $\textbf{v}_{i}, i= 1,2,\ldots,6$. If two algebras are similar, i.e. connected to each other by a transformation from the symmetry group, then their corresponding invariant solutions are connectedd to each other by the same transformation. Therfore, it is sufficient to put all similar subalgebras into one class and select a representative from each class. The set of all these representatives is called an optimal system (for details see \cite{LVO8, PO9}).\\

\section{ Symmetry group of $(3+1)-$dimensional Calogero-Bogoyavlenskii-Schiff equation}

In this section, to obtain the group transformation $\textbf{X}_i:(x,y,z,t,u)\rightarrow (\check{x}, \check{y}, \check{z}, \check{t}, \check{u})$ which is generated by the infinitesimal generator $\textbf{v}_i$ for $i=1,2,\ldots,6$; we need to solve following system of ordinary differential equations:
\begin{align*}
\frac{d(\check{x}, \check{y}, \check{z}, \check{t}, \check{u})}{d\varepsilon}= (\xi^1, \xi^2, \xi^3, \tau, \eta),\\
(\check{x}, \check{y}, \check{z}, \check{t}, \check{u})|_{\varepsilon=0}= (\xi^1, \xi^2, \xi^3, \tau, \eta).
\end{align*}
The one-parameter group $\textbf{X}_i$ spanned by $\textbf{v}_i$ is given as follows:
\begin{eqnarray}{\label{g}}
\begin{split}
\textbf{X}_1 &: (x,y,z,t,u)\rightarrow \left(x+\varepsilon (x+\lambda(t)),y,z,t+2\varepsilon t,u+\varepsilon \left(-u+\frac{\lambda^\prime(t) y}{b} +\gamma \left(\frac{b z-d y}{b},t\right)\right) \right) \\
\textbf{X}_2 &: (x,y,z,t,u)\rightarrow \left(x+\varepsilon \lambda(t),y, z,t+\varepsilon,u+\varepsilon \left(\frac{\lambda^\prime(t) y}{b} +\gamma \left(\frac{b z-d y}{b},t\right)\right) \right) \\
\textbf{X}_3 &: (x,y,z,t,u)\rightarrow \left(x+\varepsilon \lambda(t),y+\frac{\varepsilon a t}{c},z+\varepsilon t,t,u+\varepsilon \left(\frac{x}{c}+\frac{\lambda^\prime(t) y}{b} +\gamma \left(\frac{b z-d y}{b},t\right)\right) \right) \\
\textbf{X}_4 &: (x,y,z,t,u)\rightarrow \left(x+\varepsilon (-x+\lambda(t)),y+2 \varepsilon y,z+2\varepsilon z,t,u+\varepsilon \left(u+\frac{\lambda^\prime(t) y}{b} +\gamma \left(\frac{b z-d y}{b},t\right)\right) \right)\\
\textbf{X}_5 &: (x,y,z,t,u)\rightarrow \left(x+\varepsilon \lambda(t),y, z+\varepsilon,t,u+\varepsilon \left(\frac{\lambda^\prime(t) y}{b} +\gamma \left(\frac{b z-d y}{b},t\right)\right) \right)\\
\textbf{X}_6 &: (x,y,z,t,u)\rightarrow \left(x+\varepsilon \lambda(t),y+\varepsilon, z,t,u+\varepsilon \left(\frac{\lambda^\prime(t) y}{b} +\gamma \left(\frac{b z-d y}{b},t\right)\right) \right)
\end{split}
\end{eqnarray}
The entry on the right hand side gives the transformed point $\exp(x,y,z,t,u)= (\check{x}, \check{y}, \check{z}, \check{t}, \check{u})$. If $u=f(x,y,z,t)$ is known solution of Eq.(\ref{eq1}), then by using above groups $\textbf{X}_i(i=1,2,\ldots,6)$, corresponding new solutions $u_i(i=1,2,\ldots,6)$ can be obtained as follows:

\begin{eqnarray}
\begin{split}
u_1 &= f\left(x+\varepsilon (x+\lambda(t)),y,z,t+2\varepsilon t,u+\varepsilon \left(-u+\frac{\lambda^\prime(t) y}{b} +\gamma \left(\frac{b z-d y}{b},t\right)\right) \right) \\
u_2 &= f\left(x+\varepsilon \lambda(t),y, z,t+\varepsilon,u+\varepsilon \left(\frac{\lambda^\prime(t) y}{b} +\gamma \left(\frac{b z-d y}{b},t\right)\right) \right) \\
u_3 &= f\left(x+\varepsilon \lambda(t),y+\frac{\varepsilon a t}{c},z+\varepsilon t,t,u+\varepsilon \left(\frac{x}{c}+\frac{\lambda^\prime(t) y}{b} +\gamma \left(\frac{b z-d y}{b},t\right)\right) \right) \\
u_4 &= f\left(x+\varepsilon (-x+\lambda(t)),y+2 \varepsilon y,z+2\varepsilon z,t,u+\varepsilon \left(u+\frac{\lambda^\prime(t) y}{b} +\gamma \left(\frac{b z-d y}{b},t\right)\right) \right)\\
u_5 &= f\left(x+\varepsilon \lambda(t),y, z+\varepsilon,t,u+\varepsilon \left(\frac{\lambda^\prime(t) y}{b} +\gamma \left(\frac{b z-d y}{b},t\right)\right) \right)\\
u_6 &= f\left(x+\varepsilon \lambda(t),y+\varepsilon, z,t,u+\varepsilon \left(\frac{\lambda^\prime(t) y}{b} +\gamma \left(\frac{b z-d y}{b},t\right)\right) \right)
\end{split}
\end{eqnarray}
\\

\section{Symmetry reduction and exact solutions}

In this section we shall get the similarity solution for Eq.(\ref{eq1}) by solving the reduction equations which can be find with the help of similarity variables. To find the similarity variables first we solve corresponding characteristic equations which are:
\begin{equation}{\label{eq4}}
\frac{dx}{\xi^1(x,y,z,t,u)} = \frac{dy}{\xi^2(x,y,z,t,u)}=\frac{dz}{\xi^3(x,y,z,t,u)} =\frac{dt}{\tau(x,y,z,t,u)} = \frac{dv}{\eta(x,y,z,t,u)} 
\end{equation}

$\textbf{5.1. vector field }\textbf{v}_{3}$:\\
\begin{align} \label{v}
\textbf{v}_{3} =& \quad \lambda(t) \frac{\partial}{\partial x}+\frac{a t}{c} \frac{\partial}{\partial y}+t \frac{\partial}{\partial z}+\left(\frac{x}{c}+\frac{\lambda^\prime(t) y}{b} +\gamma \left(\frac{b z-d y}{b},t\right)\right) \frac{\partial}{\partial u}
\end{align}

Assume $\lambda(t)= \gamma \left(\frac{b z-d y}{b},t\right)=0$. Then associated Lagrange system is found by comprising (\ref{v}) and (\ref{eq4})
\begin{equation*}
\frac{dx}{0} = \frac{dy}{\frac{a t}{c}} =\frac{dz}{t}=\frac{dt}{0} = \frac{du}{\frac{x}{c}} 
\end{equation*}
The similarity reduction of equation $(\ref{eq1})$ in similarity form is
\begin{eqnarray}{\label{s1}}
u(x,y,t)=\frac{x y}{a t}+ f(X,Y,T) \quad
\mbox{where}\quad
X= x, \quad T=t \quad \mbox{and} \quad Y = y-\frac{a z}{c} 
\end{eqnarray} 
are the three invariants that we obtained.\\
From Eqs. $\ref{s1}$ and $(\ref{eq1})$, we get the following partial differential equation with three independent variables
\begin{equation}{\label{0}}
T f_{XT}+f_X+ A X f_{XX}+B T f_Y f_{XX}+ C T f_{XXXY}=0
\end{equation}
where $A=\frac{b}{a}, B=\frac{b c-a d}{c}, C=\frac{c-a}{c}$. The new set of infinitesimal generator for Eq. ($\ref{0}$) by applying similarity transformation method (STM) is
\begin{eqnarray*}
\xi_X &=& (a_{1}+a_{2}-a_{3}) X + a_{2} X log T +\psi(T)\\ 
\xi_Y &=& 2 a_{3} Y + a_{4}\\
\tau_T &=& 2 a_{2} T log T+2 a_{1} T\\
\eta_f &=& -a_{2} log T f -(a_{1}+a_{2}-a_{3}) f-\frac{\psi^\prime(T) T Y +A \psi(T) Y+(2 A-1)a_{2} X Y}{B T}+\frac{a_{5} X}{T}+\delta(T)
\end{eqnarray*}
where $a_{1}$, $a_{2}$, $a_{3}, a_{4}$ and $a_{5}$ are arbitrary constants and $\psi(T), \delta(T)$ are arbitrary functions.\\

\textbf{Case 1:$a_{1}\neq 0$ and else parameters and arbitrary functions are zero.}\\

 This follows the characteristic equation for $(\ref{0})$ is given by
\begin{eqnarray*}
\frac{dX}{a_{1} X}=\frac{dY}{0}=\frac{dT}{2 a_{1} T}= \frac{df}{-a_{1}f}
\end{eqnarray*}
Further, $f$ can be written as
\begin{align}\label{s2}
f(X,Y,T)= \frac{H(r,s)}{\sqrt{T}}, \quad 
 \mbox{where} \quad r = \frac{X}{\sqrt{T}} \quad \mbox{and} \quad s=Y \quad  \mbox{are similarity variables.} 
\end{align}

Further, Eq. $(\ref{0})$  can be reduced into following $(1+1)-$dimensional partial differential equation as follows: 
\begin{eqnarray}{\label{*}}
(A-1) r H_{rr}+ B H_s H_{rr}+C H_{rrrs}=0
\end{eqnarray}
Novel Lie group similarity analysis method gives the following infinitesimals when it applies on Eq.($\ref{*}$).
\begin{eqnarray*}
\xi_r &=& b_{1} r + b_{2}\\ 
\xi_s &=& \rho(s) \\
\eta_H &=& -b_{1} H - \frac{2 b_{1}(A-1)r s}{B }+\frac{b_{2} (1-A) s}{B}+b_{3} r+b_{4}+\frac{(1-A) \rho(s) r}{B}
\end{eqnarray*}
where $b_{1},b_{2}, b_{3}$ and $b_{4}$ are arbitrary constants and $\rho(s)$ is a arbitrary function of $s$.\\

\textbf{Subcase 1: $b_{1} \neq 0$, else parameters and arbitrary functions are zero.}\\

In this case, we find the similarity solution $H(r,s)$ as:
\begin{align}\label{s3}
H(r,s)= \frac{G(\zeta)}{r}+\frac{(1-A) r s}{B} \quad 
 \mbox{where} \quad \zeta=s \quad  \mbox{is a similarity variable.} 
\end{align} 
Therefore, Eq.(\ref{s3})reduced the Eq.(\ref{*}) into the following ordinary differential equation
\begin{equation}\label{*1}
2 B G(\zeta) G^\prime(\zeta)-6 C G^\prime(\zeta)=0
\end{equation}
where $B=\frac{b c-a d}{c}, C=\frac{c-a}{c}$. Hence, we get following two solutions of Eq.(\ref{*1})
\begin{align}
G(\zeta)=& \frac{3(c-a)}{(b c-a d)}\\
\mbox{or}\\ \nonumber
G(\zeta)=& -\alpha, \quad \mbox{where $\alpha$ is a new constant.}
\end{align}
By back substitution, we get following two exact solution of Eq.(\ref{eq1})
\begin{align}
u(x,y,z,t)=& \frac{(c-d)}{(b c- a d)} \frac{x y}{t}-\frac{(a-b)}{(b c- a d)} \frac{x z}{t}+\frac{3(c-a)}{(b c-a d) x}\\
u(x,y,z,t)=& \frac{(c-d)}{(b c- a d)} \frac{x y}{t}-\frac{(a-b)}{(b c- a d)} \frac{x z}{t}+\frac{\alpha}{x} .
\end{align}
\\
\textbf{Subcase 2: $b_{2} \neq 0, \rho(s) \neq 0$, else parameters are zero.}\\

In this case, Lagrange's characteristic equations are as 
\begin{eqnarray*}
\frac{dr}{b_{2}}=\frac{ds}{\rho(s)} = \frac{dH}{\frac{b_{2} (1-A) s}{B}+\frac{(1-A) \rho(s) r}{B}}
\end{eqnarray*}
 Hence, we find the similarity solution $H(r,s)$ as:
\begin{align}
H(r,s)&= G(\zeta)+\frac{(1-A) s r}{B},   \quad 
 \mbox{where}\\ \nonumber
  \zeta &= r- b_{2} \int \frac{ds}{\rho(s)} \quad  \mbox{is a similarity variable.} 
\end{align} 
Therefore, Eq.(\ref{*}) has reduced into the following ordinary differential equation
\begin{equation}\label{*2}
G^{\prime \prime \prime \prime}(\zeta)=0
\end{equation}
Hence, solution of Eq.(\ref{*2}) is
\begin{align}
G(\zeta)= \beta_1+ \beta_2 \zeta+\beta_3 \zeta^2+\beta_4 \zeta^3
\end{align}
By back substitution, exact solution of Eq.(\ref{eq1})is
\begin{align}
u(x,y,z,t)=& \frac{(c-d)}{(b c- a d)} \frac{x y}{t}-\frac{(a-b)}{(b c- a d)} \frac{x z}{t}+\frac{\beta_1 \sqrt{t}+ \beta_2 x}{(t}-\frac{b_2 \beta_2}{\sqrt{t}} \int \frac{dY}{\rho(Y)}\\ \nonumber 
+&\frac{ \beta_3}{\sqrt{t}} \left[\frac{x}{\sqrt{t}} - \int \frac{ b_2 dY}{\rho(Y)}\right]^2+\frac{ \beta_4}{\sqrt{t}} \left[\frac{x}{\sqrt{t}} - \int \frac{ b_2 dY}{\rho(Y)}\right]^3\\ \nonumber
\mbox{here} \quad Y= y- \frac{a z}{c}
\end{align}
\\
$\textbf{5.2. vector field }\textbf{v}_{1}$:\\
\begin{align}\label{vm}
\textbf{v}_{2} =& \quad  \lambda(t) \frac{\partial}{\partial x}+\frac{\partial}{\partial t}+\left(\frac{\lambda^\prime(t) y}{b} +\gamma \left(\frac{b z-d y}{b},t\right)\right)\frac{\partial}{\partial u}
\end{align}
Assume $ \gamma \left(\frac{b z-d y}{b},t\right)=\frac{b z-d y}{b}$. Then corresponding Lagrange's characteristic equations can be find by comprising ($\ref{vm}$) and (\ref{eq4})
\begin{equation}
\frac{dx}{\lambda(t)} = \frac{dy}{0} =\frac{dz}{0}=\frac{dt}{1} = \frac{du}{\frac{\lambda^\prime(t) y}{b} +\frac{b z-d y}{b}} 
\end{equation}
The similarity reduction of equation $(\ref{eq1})$ in similarity form with following three similarity variables is
\begin{eqnarray}{\label{s21}}
u(x,y,t)&=& f(X,Y,Z)+\frac{\lambda(t) y}{b} +\frac{(b z-d y) t}{b}\\ \nonumber
X &=& x-\int\lambda(t)dt \quad Z=z \quad \mbox{and} \quad Y = y .
\end{eqnarray} 

From Eqs. $\ref{s21}$ and $(\ref{eq1})$, we get the following partial differential equation with three independent variables
\begin{equation}{\label{00}}
a f_X f_{XY}+a f_Y f_{XX}+c f_X f_{XZ}+d\, f_Z f_{XX}+ f_{XXXY}+f_{XXXZ}=0
\end{equation}
 The new set of infinitesimal generator for Eq. ($\ref{00}$) by applying similarity transformation method (STM) is
\begin{align*}
\xi_X &= a_{4} X + a_{5}\\ 
\xi_Y &=  a_{1} Y + a_{3}\\
\tau_T &=  a_{1} Z+ a_{1} \\
\eta_f &=-a_{4} f + \psi\left(\frac{b Z-d Y}{b}\right) 
\end{align*}
where $a_{1}$, $a_{2}$, $a_{3}, a_{4}$ and $a_{5}$ are arbitrary constants and $\psi(\frac{b Z-d Y}{b})$ is an arbitrary function.\\

Assume $\psi(\frac{b Z-d Y}{b})=0$. Then Lagrange's characteristic equation for $(\ref{00})$ is given by
\begin{eqnarray*}
\frac{dX}{a_{4} X + a_{5}}=\frac{dY}{ a_{1} Y + a_{3}}=\frac{dT}{ a_{1} Z+ a_{1}}= \frac{df}{-a_{4}f}
\end{eqnarray*}
Further, $f$ can be written as
\begin{align}\label{s22}
f(X,Y,Z)= \frac{H(r,s)}{(a_{1} Y+a_{3})^{\frac{a_{4}}{a_{1}}}} \, ,\quad  
 \mbox{where} \quad r = \frac{a_{4} X+a_{5}}{(a_{1} Y+a_{3})^{\frac{a_{4}}{a_{1}}}} \quad \mbox{and} \quad s=\frac{a_{1} Z+a_{2}}{(a_{1} Y+a_{3})} \quad  \mbox{are similarity variables.} 
\end{align}

Further, Eq. $(\ref{00})$  can be reduced into following $(1+1)-$dimensional partial differential equation as follows: 
\begin{eqnarray}\label{p}
&-(a+b) a H_{r}H_{rr}+(c-a s) \frac{a_{1}}{a_{4}} H_{r} H_{rs}-b H H_{rr} +(d-b s)\frac{a_{1}}{a_{4}} H_{s} H_{rr}\\ \nonumber &- 2 a H_{r}^2-4 a_{4} H_{rrr}-a_{4} r H_{rrrr} + a_{1} (1-s) H_{rrrs}=0
\end{eqnarray}
Novel Lie group analysis method gives the following infinitesimals when it applies on Eq.($\ref{p}$).
\begin{align*}
\xi_r &= b_{1} r \\ 
\xi_s &= 0 \\
\eta_H &= -b_{1} H - b_{2} (b s-d)^{\frac{-a_{4}}{a_{1}}}
\end{align*}
where $b_{1}$ and $b_{2}$ are arbitrary constants.\\
 The similarity solution $H(r,s)$ can be written in the following similarity form:
\begin{align}
H(r,s) = \frac{G(\zeta)}{r}+\frac{b_{2} (b s-d)^{\frac{-a_{4}}{a_{1}}}}{b_{1}} \, , \quad 
\mbox{with similarity variable} \quad
 \zeta =s .
\end{align}
This transformation reduces Eq.(\ref{eq1}) into following ordinary differential equation
\begin{equation} \label{ss}
(c+d-(a+b)s) G(\zeta) G^\prime(\zeta)+6 a_{4} (1-s)G^\prime(\zeta)=0
\end{equation}
Therefore, we get following two solution of Eq.(\ref{ss})
\begin{align}
G(\zeta)&= \frac{6 a_{4} (1-\zeta)}{(a+b) \zeta-(c+d)}\\
G(\zeta) &= \beta
\end{align}
Hence, we find the solutions of $(\ref{eq1})$ which are given by
\begin{align}
u(x,y,z,t)&= \frac{\lambda(t)y}{b}+\frac{(b z-d y) t}{b}+\frac{b_2\left(b(\frac{a_1 z+a_2}{a_1 y+a_3})-d\right)^{\frac{-a_4}{a_1}}}{a_1 (a_1 y+a_3)^{\frac{a_4}{a_1}}} + \frac{6 a_4(1-\frac{a_1 z+a_2}{a_1 y+a_3})}{[a_5+a_4 x-a_4 \int \lambda(t) dt][(a+b)\frac{a_1 z+a_2}{a_1 y+a_3}-(c+d)]}  \\
u(x,y,z,t)&= \frac{\lambda(t)y}{b}+\frac{(b z-d y) t}{b}+\frac{b_2\left(b(\frac{a_1 z+a_2}{a_1 y+a_3})-d\right)^{\frac{-a_4}{a_1}}}{a_1 (a_1 y+a_3)^{\frac{a_4}{a_1}}}+\frac{\beta}{a_5+a_4 x-a_4 \int \lambda(t) dt}
\end{align}
\section{ Conclusion}
 In this paper we have shown that the Calogero-Bogoyavlenskii-Schiff equation can be transformed by a Lie point transformation to new partial differential equation with less independent variables and again Lie group symmetry reduces these equation into ordinary differential equation. Infinitesimal generators  and different vector fields for CBS equation  were obtained by using the Lie symmetry group analysis. Using the criterion of invariance of Eq.(\ref{eq1}) under the the infinitesimal prolongation, we find the Lie symmetries group of $(3+1)$-dimensional Calogero-Bogoyavlenskii-Schiff equation and similarity variables played a important role in the reduction of equation. This work is significant since the exact solutions  obtained shall be helpful  in other applied sciences as condensed matter physics, field theory, fluid dynamics, plasma physics, non-linear optics, etc. Our exact solutions may serve as benchmark in the accuracy testing and comparison of their numerical algorithms. The availability of computer systems like Mathematica or Maple facilitates the tedious algebraic calculations. The method which we have proposed in this article is also a standard, direct and computer-literate method, which allows us to solve complicated and tedious algebraic calculation.\\

Acknowledgments:  The first author also expresses her gratitude to the University Grants Commission, New Delhi, India for financial support to carry out the above work.

\end{document}